\documentclass[12pt]{iopart}

\usepackage{iopams}  
\usepackage{subfigure}

\usepackage{graphicx}

\begin{document}

\title{Approximate matrix completion based on cavity method}

\author{Chihiro Noguchi and Yoshiyuki Kabashima}

\address{Department of Mathematical and Computing Science,
Tokyo Institute of Technology,
2-12-1, Ookayama, Meguro-ku, Tokyo, Japan}
\ead{noguchi.c.aa@m.titech.ac.jp}
\begin{abstract}
In order to solve large matrix completion problems with practical computational cost, an approximate approach based on matrix factorization has been widely used. Alternating least squares (ALS) and stochastic gradient descent (SGD) are two major algorithms to this end. In this study, we propose a new algorithm, namely cavity-based matrix factorization (CBMF) and approximate cavity-based matrix factorization (ACBMF), which are developed based on the cavity method from statistical mechanics. ALS yields solutions with less iterations when compared to those of SGD. This is because its update rules are described in a closed form although it entails higher computational cost. CBMF can also write its update rules in a closed form, and its computational cost is lower than that of ALS. ACBMF is proposed to compensate a disadvantage of CBMF in terms of relatively high memory cost. We experimentally illustrate that the proposed methods outperform the two existing algorithms in terms of convergence speed per iteration, and it can work under the condition where observed entries are relatively fewer. Additionally, in contrast to SGD, (A)CBMF does not require scheduling of the learning rate.

\end{abstract}

\maketitle


\section{Introduction}

Recent technological advances triggered the generation and accumulation of significant amounts of data. In response to the trend, several methods are proposed to extract useful information from them. This produced significant results in various fields including science and engineering. A typical example can be found in collaborative filtering, which is a methodology that is used in recommender systems \cite{koren2009matrix}. As a comprehensive example, we consider a user-movie matrix $Y\in\mathbb{R}^{N\times M}$, where $N$ and $M$ denote the number of users and movies, respectively, and an entry of $Y$, $Y_{ij}$, denotes rating from user $i$  movie $j$. Users normally evaluate only a small fraction of movies, and thus most entries of $Y$ are missing. Under the aforementioned types of setting, the primary objective of matrix completion involves predicting  missing entries.

A natural approach for this involves minimizing the rank of the matrix under constraints yielded by observed entries, and this is generally referred to as ``low-rank matrix completion''. Unfortunately, it is NP-hard to literally solve the rank minimization problem. In order to practically overcome the difficulty, relaxation of matrix rank to nuclear norm was proposed \cite{fazel2002matrix}. Interestingly, it is guaranteed that the solution of the nuclear norm minimization is exactly in agreement with that of the original rank minimization if certain conditions are satisfied \cite{candes2009exact,candes2010power,recht2011simpler,keshavan2010matrix}. The minimization of nuclear norm belongs to the class of convex optimization problems, and thus the optimal solution is determined via versatile semidefinite programming solvers when the matrix size is relatively small. However, in several realistic problems, matrix sizes are not so small, and computational and memory costs required by the nuclear norm minimization often exceed practically acceptable levels.

In order to deal with such situations, a non-convex approach using matrix factorization was proposed more recently \cite{koren2009matrix}. When the objective matrix is factorized into two matrices of lower rank, nuclear norm is evaluated as the sum of their Frobenius norms. The non-convex formulation significantly reduces necessary computational and memory costs while we can generally find only local minima. However, a recent study \cite{ge2016matrix} indicated that under a certain condition, the objective function of matrix factorization does not exhibit spurious local minima. Each local minimum is transformed to another via trivial operations such as permutations of column/rows with high probabilities. 

Two major algorithms, alternating least squares (ALS) \cite{zhou2008large,jain2013low,hardt2014understanding} and stochastic gradient descent (SGD) \cite{yu2012scalable,zhuang2013fast,recht2013parallel} are proposed for the matrix factorization to date. The main objective of this study is to develop a new algorithm by borrowing an idea from the cavity method from statistical mechanics. Even if the absence of spurious local minima is guaranteed, the performance of the solution search is determined via dynamical properties of the used algorithm. We experimentally illustrate that the proposed cavity-based algorithms exhibit better performance than the two algorithms without delicate tuning of control parameters when the number of observed data is relatively small. 

Several extant studies apply the cavity method for the matrix factorization problems. An approximate message passing (AMP) based approach to generalized bilinear inference problem including the matrix completion was proposed in \cite{parker2014bilinear,parker2014bilinear2}. A detailed derivation of AMP-type algorithms and performance analysis for the Bayes optimal cases are provided in \cite{kabashima2016phase}. Reference \cite{matsushita2013low} presents an AMP based algorithm for low-rank matrix reconstruction and its application to {\it K}-means type clustering. All of these methods follow the Bayesian framework.
The differences of the present study from these are as follows. We do not employ the Bayesian approach, and thus it is not necessary to select a prior distribution. Additionally, we focus on the matrix completion as a particular application of matrix factorization, and aim to develop efficient algorithms exploiting the properties of the specific problem.

The remainder is organized as follows. In section 2, the problem setting is detailed. In section 3, we explain the details of the proposed algorithm. In section 4, the performance of the proposed algorithms is illustrated via applications for synthetic and realistic data. The final section presents the summary.

\section{Problem Setting}

In the simplest case, low rank matrix completion is defined as follows:

\begin{eqnarray}
&\min \limits_{X}&\ \ \ \ \ \ \mathrm{rank}(X) \nonumber\\
&\mathrm{subject\ to}&\ \ \ \ \ \  X_{\mu i} = Y_{\mu i}\ \ \ (\mu,i)\in\Omega\ ,
\label{eq:rank_minimization}
\end{eqnarray}

\noindent
where $X$ and $Y$ are denoted as decision variables and observed entries, respectively, and $\Omega$ stands for the set of indices of $Y$. The problem is guaranteed to exhibit a unique solution with high probability when the size of $\Omega$ is sufficiently large. However, there is no known algorithm that solves (\ref{eq:rank_minimization}) in practical time. Hence, relaxation of the matrix rank to the nuclear norm is typically employed and defined as follows:

\begin{equation}
\|X\|_*=\sum_k^{\min\{N,M\}} \sigma_k\ ,
\end{equation}

\noindent
where $\sigma_k$ denotes the $k$th highest singular value of $X$. In Lagrange form, the nuclear norm relaxation converts (\ref{eq:rank_minimization}) as follows:

\begin{equation}
\label{eq:nuclear_noise}
\min\limits_{X\in\mathbb{R}^{N\times M}}\frac{1}{2}\sum_{(\mu,i)\in\Omega}(Y_{\mu i}-X_{\mu i})^2 + \lambda\|X\|_*\ .
\end{equation}

\noindent
The solution of (\ref{eq:nuclear_noise}) is determined in a polynomial time via versatile solvers of semi-definite programming. However, such solvers require singular value decomposition per iteration, and their computational and memory costs easily exceed practically permissible levels when the system size increases.

A popular approach to overcome this disadvantage involves using non-convex relaxation. Let us assume that rank of $X$ is $R$, and this means that $X$ is expressed as $X=UV^T$ by using two smaller matrices as $U\in\mathbb{R}^{N\times R},V\in\mathbb{R}^{M\times R}$. An attractive property of the nuclear norm is that it is evaluated by another norm as follows:

\begin{equation}
\label{eq:nuclear_factorization}
\|X\|_* = \inf \left\{ \frac{1}{2}\|U\|_F + \frac{1}{2}\|V\|_F\ :\ X=UV^T \right\}\ ,
\end{equation}

\noindent
where $\|A\|_F=\sqrt{\sum_{ij}A_{ij}^2}$ denotes the Frobenius norm of matrix $A$ \cite{recht2010guaranteed}. We insert (\ref{eq:nuclear_factorization}) into (\ref{eq:nuclear_noise}) to yield a non-convex version of (\ref{eq:nuclear_noise}) as

\begin{equation}
\label{eq:nonconvex_relaxed}
\min\limits_{U\in\mathbb{R}^{N\times R},V\in\mathbb{R}^{M\times R}}\frac{1}{2}\sum_{(\mu,i)\in\Omega}\left(Y_{\mu i}-\sum_{r=1}^R  U_{\mu r}V_{ir}\right)^2 + \frac{1}{2}\lambda\|U\|_F^2 + \frac{1}{2}\lambda\|V\|_F^2\ .
\end{equation}

\noindent
In contrast to (\ref{eq:nuclear_noise}), (\ref{eq:nonconvex_relaxed}) ceases to be convex, and this implies that multiple local minima can exist. However, it was recently illustrated that a spurious local minimum is absent with a high probability if a few conditions are satisfied \cite{ge2016matrix}.

Without any constraints, the degree of freedom of this problem is given as $R(N+M)$. This implies that the number of observations $|\Omega|$ must not be less than $R(N+M)$ to determine a solution. In the following, we assume that this condition is satisfied.

Two major algorithms, ALS and SGD, are known to solve (\ref{eq:nonconvex_relaxed}). ALS is widely known as a standard approach to non-convex optimization problems due to its simplicity. When $V$ is fixed, each row of $U$ is independently calculated, and the objective function (\ref{eq:nonconvex_relaxed}) is then expressed as follows:

\begin{equation}
\min\limits_{{\bf u}_\mu} \frac{1}{2} \sum_{i\in\partial\mu}(y_{\mu i}-{\bf u}_\mu^T{\bf v}_i)^2+\lambda\|{\bf u}_\mu\|^2 ,
\label{eq:als_quad}
\end{equation}

\noindent
where ${\bf u}_\mu$ and ${\bf v}_i$ denote the $\mu$-th and $i$-th rows of $U$ and $V$ respectively, and $\partial\mu$ denotes a set of observed indices of $\mu$-th row of $Y$. Thus, (\ref{eq:als_quad}) leads to the following closed form solution:

\begin{equation}
{\bf u}_\mu^*=\left( \sum_{i\in\partial\mu} {\bf v}_i{\bf v}_i^T+\lambda{\bf I}_R \right)^{-1}\left(\sum_{i\in\partial\mu}y_{\mu i}{\bf v}_i\right) , \label{eq:als_inv}
\end{equation}

\noindent
where ${\bf I}_R $ denotes $R\times R$ unit matrix. Subsequently, we fix $U$ and solve $V$ in turn, and ALS repeats this operation until convergence. The main advantage of ALS is the ease of parallelization although the computational cost per iteration exceeds that of SGD.

The other algorithm, SGD, is also widely known as a standard algorithm for continuous optimization problems. Specifically, SGD computes a gradient only with respect to pairwise indices $(\mu, i)\in\Omega$ selected at random per iteration, and the gradient updates the corresponding variables based on the given learning rate $\eta$ as follows:

\begin{eqnarray}
{\bf u}_\mu&\leftarrow&{\bf u}_\mu-\eta\left\{\lambda{\bf u}_\mu-(y_{\mu i}-{\bf u}_\mu{\bf v}_i)^T{\bf v}_i\right\} , \\
{\bf v}_i&\leftarrow&{\bf v}_i-\eta\left\{\lambda{\bf v}_i-(y_{\mu i}-{\bf u}_\mu{\bf v}_i)^T{\bf u}_\mu\right\} ,
\end{eqnarray}

\noindent
The algorithm exhibits an advantage wherein its computational cost in the elemental update is lower. However, it has two major disadvantages. The first is that an overwriting issue can arise when the several updates are conducted in parallel. The second is that it is highly sensitive to the learning rate. Distributed SGD (DSGD) \cite{recht2013parallel} (the name {\it Jellyfish} used in \cite{recht2013parallel}) overcomes the first disadvantage by dividing the observed matrix into a few blocks, considering a set of independent blocks, and updating a pair of indices from each block in it. However, the second disadvantage still remains, and the learning rate should be carefully tuned and scheduled. The adjustment of the learning rate significantly affects the convergence of the algorithm.

\begin{figure}[htbp]
  \begin{center}
    \begin{tabular}{c}

      \begin{minipage}{0.57\hsize}
        \begin{center}
          \includegraphics[clip, width=9cm]{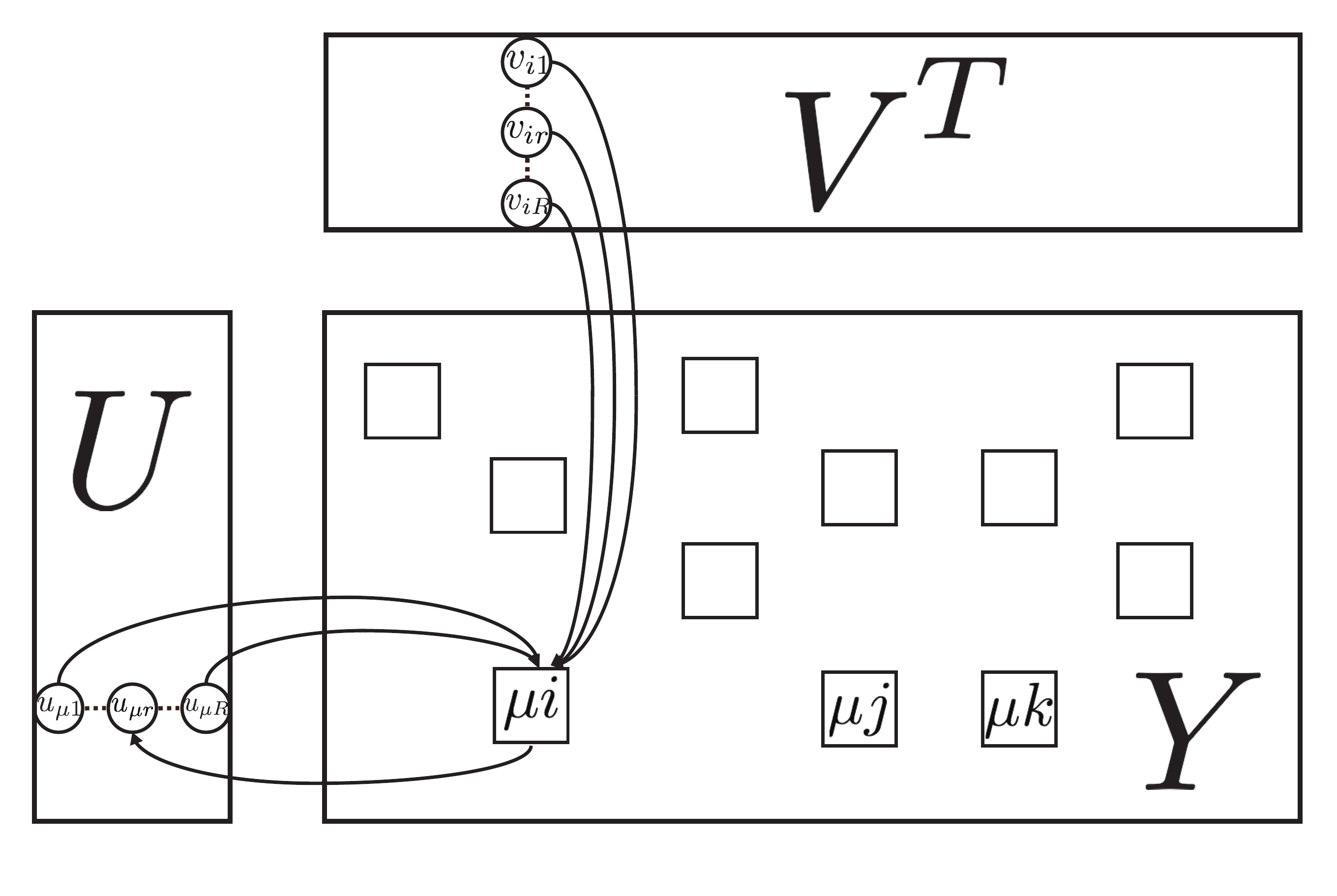}
          \hspace{1.6cm}
        \end{center}
      \end{minipage}
    \hspace{-20pt}
      \begin{minipage}{0.57\hsize}
        \begin{center}
          \includegraphics[clip, width=5.5cm]{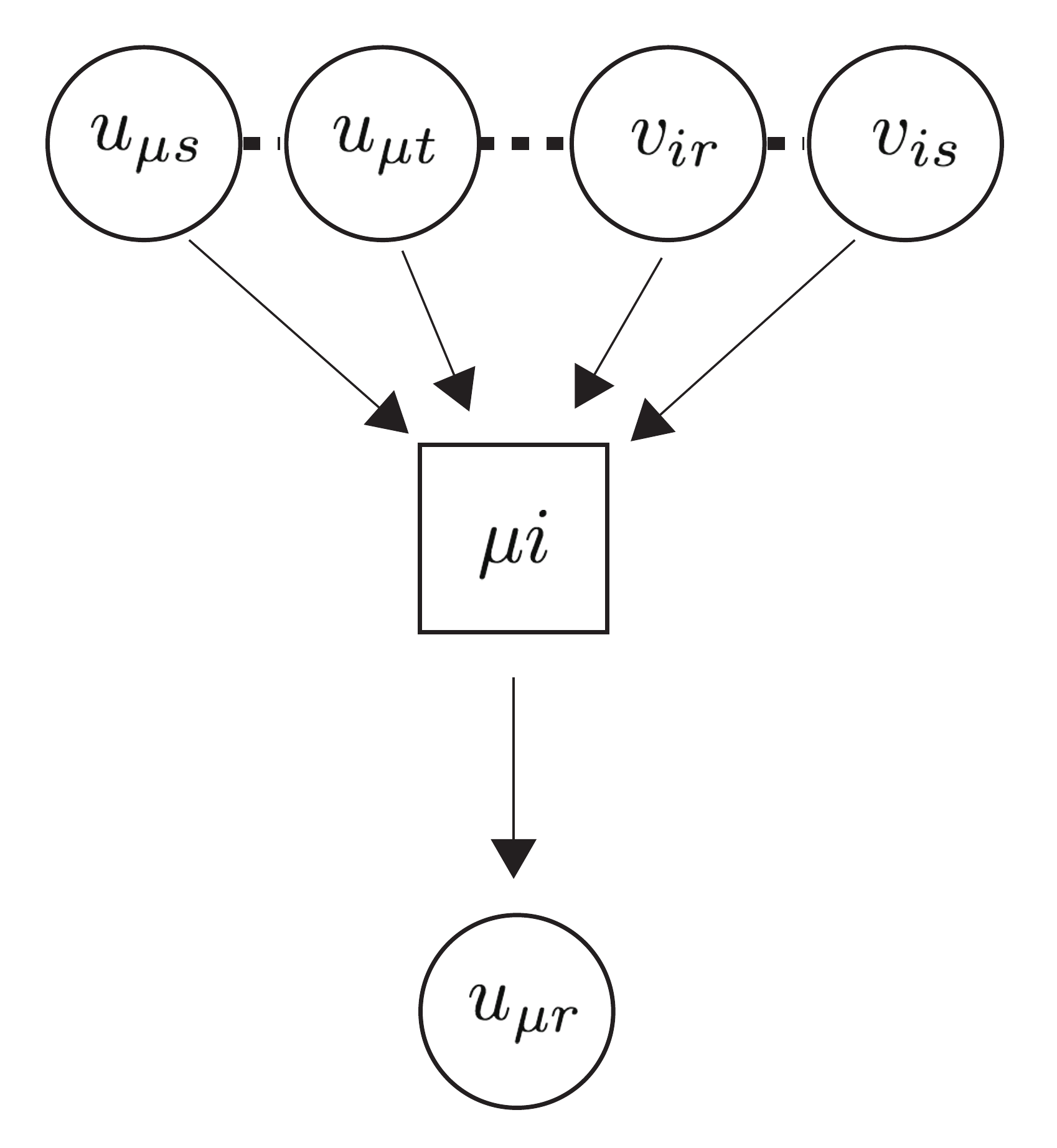}
          \hspace{1.6cm}
        \end{center}
      \end{minipage}

    \end{tabular}
    \caption{Graphical expression of (\ref{eq:bp1}) (left) and its enlarged illustration (right). Circle and squares correspond to variable and function nodes, respectively. Equation (\ref{eq:bp2}) is also computed in a similar manner.}
    \label{fig:ftov}
  \end{center}
\end{figure}

\begin{figure}[htbp]
  \begin{center}
    \begin{tabular}{c}

      \begin{minipage}{0.57\hsize}
        \begin{center}
          \includegraphics[clip, width=9cm]{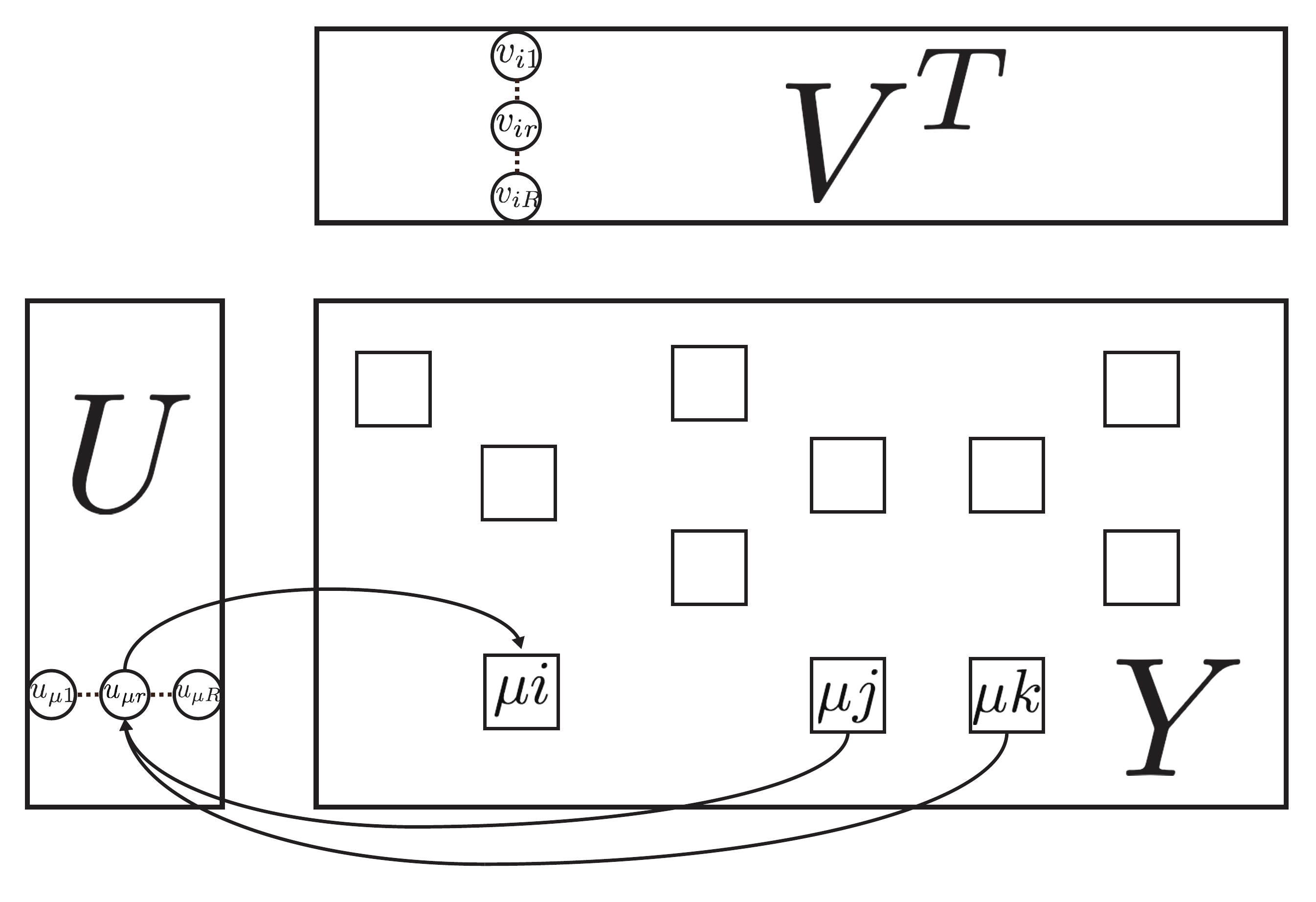}
          \hspace{1.6cm}
        \end{center}
      \end{minipage}
      
    \hspace{-20pt}
      \begin{minipage}{0.57\hsize}
        \begin{center}
          \includegraphics[clip, width=5cm]{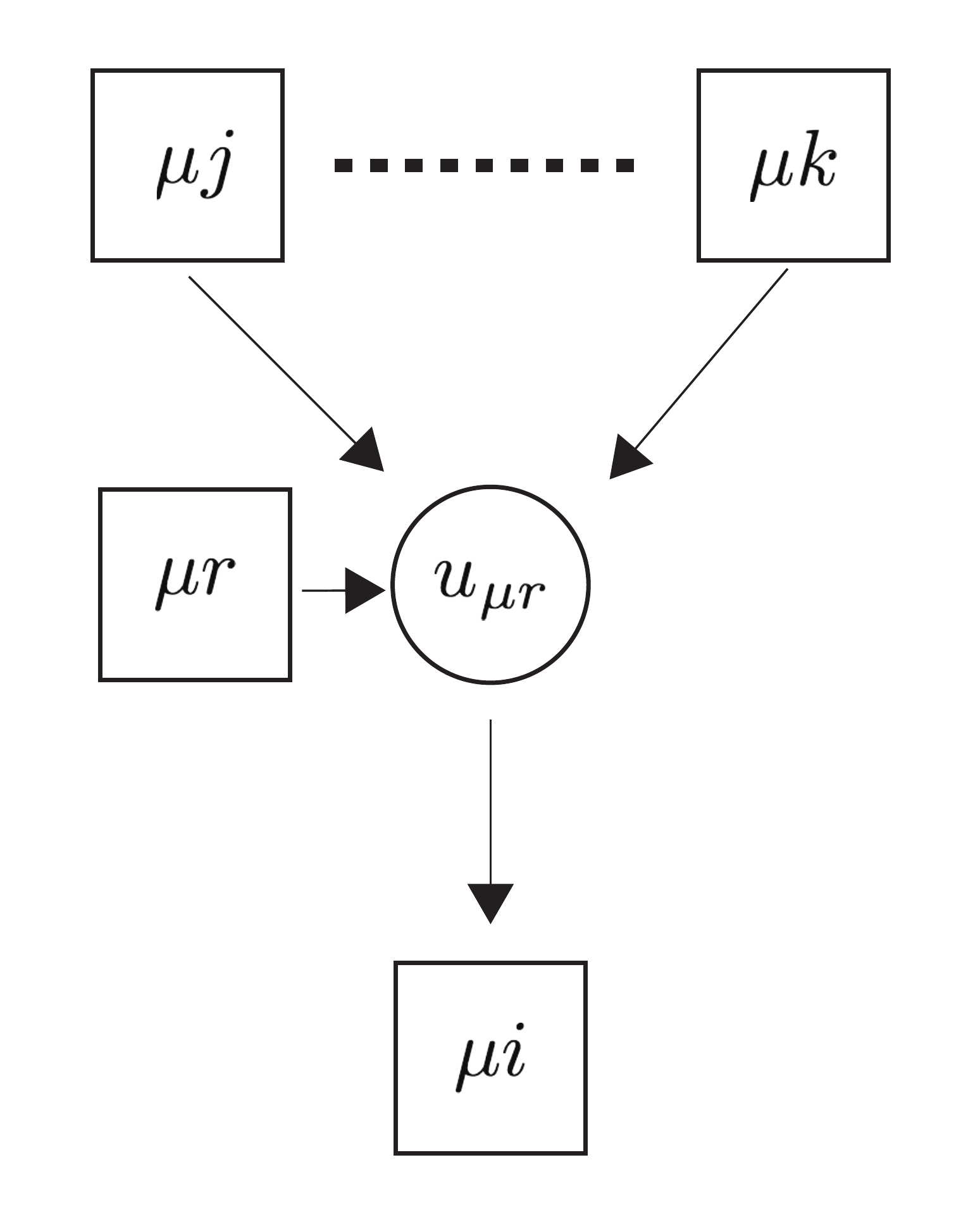}
          \hspace{1.6cm}
        \end{center}
      \end{minipage}

    \end{tabular}
    \caption{Graphical expression of (\ref{eq:bp3}) (left) and its enlarged illustration (right). Equation (\ref{eq:bp4}) is also computed in a similar manner.}
    \label{fig:vtof}
  \end{center}
\end{figure}

\section{A Cavity-Based Approach}

In order to explore the possibility of achieving a better performance, we develop an algorithm for the matrix factorization based on the cavity method \cite{mezard2009information}. Thus, we first express the variable dependence of (\ref{eq:nonconvex_relaxed}) by a factor graph (Figure \ref{fig:ftov}). The variable nodes are expressed by circles and denote entries of two matrices $U$ and $V$ while the factor nodes are represented by squares and stand for factors constituting (\ref{eq:nonconvex_relaxed}), namely, $(1/2)\left( Y_{\mu i}-\sum_{r=1}^R U_{\mu r}V_{ir} \right)^2$,$(\lambda/2)U_{\mu r}^2$ and $(\lambda/2)V_{ir}^2$. An edge for a pair of variable and factor nodes is provided if and only if the variable and factor nodes are directly related.

The basic idea of the cavity method is to approximate the multivariate minimization problem (\ref{eq:nonconvex_relaxed}) via a bunch of minimization problems with respect to single variables.
Hence, we introduce ``cavity objective functions'' $f_{\mu r\rightarrow(\mu i)}(u_{\mu r})$ and $g_{ir\rightarrow(\mu i)}(v_{ir})$. 
The function $f_{\mu r\rightarrow(\mu i)}(u_{\mu r})$ denotes the objective function after the minimization with respect to all variables other than $u_{\mu r}$ is performed in the ``$(\mu i)$-cavity system'' that is defined by removing $(1/2)\left(Y_{\mu i}-\sum_{r=1}^R U_{\mu r}V_{ir}\right)^2$ from (\ref{eq:nonconvex_relaxed}), and similarly for $g_{ir\rightarrow(\mu i)}(v_{ir})$. 
The summation of the cavity objective functions and $(1/2)\left(Y_{\mu i}-\sum_{r=1}^R U_{\mu r}V_{ir}\right)^2$ approximates the full objective function of (\ref{eq:nonconvex_relaxed}).
Conversely, we remove the contribution of  $f_{\mu r\rightarrow(\mu i)}(u_{\mu r})$ from the full summation and minimize the resulting function with respect to all variables except for $u_{\mu r}$. 
This yields ``cavity bias function'' $\hat{f}_{(\mu i)\rightarrow\mu r}(u_{\mu r})$, and this denotes the effective influence of the factor $(1/2)\left(Y_{\mu i}-\sum_{r=1}^R U_{\mu r}V_{ir}\right)^2$ to the variable $u_{\mu r}$, and similarly for $\hat{g}_{(\mu i)\rightarrow ir}(v_{ir})$. 
The summation of the cavity bias functions with the exception of $\hat{f}_{(\mu i)\rightarrow\mu r}(u_{\mu r})$ and $(\lambda/2)u_{\mu r}^2$ yields $f_{\mu r\rightarrow(\mu i)}(u_{\mu r})$, and similarly for $g_{ir\rightarrow(\mu i)}(v_{ir})$. 
They constitute a closed set of functional equations to determine the cavity objective and bias functions as follows:

\begin{equation}
\fl
\hat{f}_{(\mu i)\rightarrow\mu r}(u_{\mu r})
=
\min\limits_{\{{\bf u}_{\mu},{\bf v}_{i}\}\backslash u_{\mu r}}
\left\{
\frac{1}{2}(y_{\mu i}-\sum_s u_{\mu s}v_{is})^2
+\sum_{s\neq r} f_{\mu s\rightarrow (\mu i)}(u_{\mu s})
+\sum_{s} g_{is\rightarrow (\mu i)}(v_{is})
\right\} , \label{eq:bp1}
\end{equation}

\begin{equation}
\fl
\hat{g}_{(\mu i)\rightarrow i r}(v_{i r})
=
\min\limits_{\{{\bf u}_{\mu},{\bf v}_{i}\}\backslash v_{i r}}
\left\{
\frac{1}{2}(y_{\mu i}-\sum_s u_{\mu s}v_{is})^2
+\sum_{s} f_{\mu s\rightarrow (\mu i)}(u_{\mu s})
+\sum_{s\neq r} g_{is\rightarrow (\mu i)}(v_{is})
\right\} , \label{eq:bp2}
\end{equation}

\begin{equation}
f_{\mu r\rightarrow(\mu i)}(u_{\mu r})
=
\sum_{(\mu j)\in\partial\mu r\backslash (\mu i)} \hat{f}_{(\mu j)\rightarrow \mu r}(u_{\mu r})
+ \frac{1}{2}\lambda u_{\mu r}^2 , \label{eq:bp3}
\end{equation}

\begin{equation}
g_{i r\rightarrow(\mu i)}(v_{i r})
=
\sum_{(\nu i)\in\partial i r\backslash (\mu i)} \hat{g}_{(\nu i)\rightarrow i r}(v_{i r})
+ \frac{1}{2}\lambda v_{i r}^2 , \label{eq:bp4}
\end{equation}

\noindent
where ${\bf u}_\mu$ and ${\bf v}_i$ denote the $\mu$-th and $i$-th rows of $U$ and $V$, respectively, and $A\backslash a$ generally indicates a set that is defined via eliminating an element $a$ from a set $A$. The indices of factor nodes are denoted with parentheses while those of variable nodes are not. The notation $\partial\mu r$ stands for the set of factor nodes that directly connect variable node indexed by $\mu r$. After determining the cavity objective and bias functions from (\ref{eq:bp1})-(\ref{eq:bp4}), ``marginal'' objective functions for each variable are provided as follows:

\begin{equation}
f_{\mu}(u_{\mu r})
=
\sum_{(\mu i)\in\partial\mu r} \hat{f}_{(\mu i)\rightarrow \mu r}(u_{\mu r})
+ \frac{1}{2}\lambda u_{\mu r}^2 , \label{eq:bp_marge1}
\end{equation}

\begin{equation}
g_{i}(v_{i r})
=
\sum_{(\mu i)\in\partial i r} \hat{g}_{(\mu i)\rightarrow i r}(v_{i r})
+ \frac{1}{2}\lambda v_{i r}^2 . \label{eq:bp_marge2}
\end{equation}

\noindent
Thus, entries of the factorized matrices are evaluated as follows:

\begin{equation}
u_{\mu r}^*
=
\mathop{\rm arg~min}\limits_{u_{\mu r}} \left\{f_{\mu}(u_{\mu r})\right\} , \label{eq:fixedu}
\end{equation}

\begin{equation}
v_{ir}^*
=
\mathop{\rm arg~min}\limits_{v_{ir}} \left\{g_{i}(v_{ir})\right\} . \label{eq:fixedv}
\end{equation}

\subsection{Derivation of the algorithm}
Two issues are emphasized here. First, when the factor graph does not contain any cycles, the solution given by the cavity method is exact. However, cycles generally exist in the matrix factorization problem. However, if the positions of the observed entries are randomly selected and their number is limited up to $O(N)$ as assumed in the following, then the resulting factor graph is considered as a sparse random graph. Thus, the lengths of the cycles typically scale as $O(\ln N)$ when the system size $N$ increases. Therefore, it is reasonable to expect that the cavity method yields reasonably accurate approximates for large $N$ as the effect of the cycles becomes negligible. Second, solving (\ref{eq:bp1})-(\ref{eq:bp4}) is, unfortunately, technically difficult since they are provided as functional equations. In order to overcome the difficulty, we parameterize the cavity objective and bias functions in the form of quadratic functions as follows:

\begin{eqnarray}
\hat{f}_{(\mu i)\rightarrow\mu r}(u_{\mu r})&=&\frac{1}{2}\hat{a}_{(\mu i)\rightarrow\mu r}u_{\mu r}^2-\hat{b}_{(\mu i)\rightarrow\mu r}u_{\mu r} , \label{eq:equad1} \\
\hat{g}_{(\mu i)\rightarrow i r}(v_{i r})&=&\frac{1}{2}\hat{c}_{(\mu i)\rightarrow ir}v_{ir}^2-\hat{d}_{(\mu i)\rightarrow ir}v_{ir} , \label{eq:equad2}\\
f_{\mu r\rightarrow(\mu i)}(u_{\mu r}) &=&\frac{1}{2}a_{\mu r\rightarrow(\mu i)}u_{\mu r}^2-b_{\mu r\rightarrow(\mu i)}u_{\mu r} + \frac{1}{2}\lambda u_{\mu r}^2 ,\label{eq:equad3} \\
g_{ir\rightarrow(\mu i)}(v_{ir}) &=&\frac{1}{2}c_{ir\rightarrow(\mu i)}v_{ir}^2-d_{ir\rightarrow(\mu i)}v_{ir} + \frac{1}{2}\lambda v_{ir}^2 . \label{eq:equad4}
\end{eqnarray}

\noindent
However,  the insertion of (\ref{eq:equad1})-(\ref{eq:equad4}) into (\ref{eq:bp1})-(\ref{eq:bp4}) does not yield a closed form of equations to determine the parameters. This indicates that a further approximation is required. Hence, we assign ${\bf v}_i$ the value in the previous step to solve the minimization problem of (\ref{eq:bp1}). Similarly for equation (\ref{eq:bp2}). This leads to quadratic forms with respect to ${\bf u}^r_\mu$ and ${\bf v}^r_i$ from (\ref{eq:bp1}) and (\ref{eq:bp2}), respectively. Here, ${\bf u}^r_\mu$  denotes a vector excluding $u_{\mu r}$ from ${\bf u}_\mu$. Similarly, this stands for ${\bf v}^r_i$ . Accordingly, when ${\bf v}_i$ is fixed, the equation $(\ref{eq:bp1})$ is re-expressed as follows:

\begin{eqnarray}
\fl \nonumber
\hat{f}_{(\mu i)\rightarrow\mu r}(u_{\mu r}) = \\\fl
\min\limits_{\{{\bf u}_{\mu}\}\backslash u_{\mu r}}
\left\{
\frac{1}{2}({\bf u}^r_\mu)^T\left(\Gamma_{{\bf a}^r_{\mu\rightarrow(\mu i)}} + {\bf v}^r_i({\bf v}^r_i)^T\right){\bf u}^r_\mu - \left\{ {\bf b}_{\mu\rightarrow(\mu i)}^r + (y_{\mu i}-u_{\mu r}v_{ir}){\bf v}_i\right\}^T {\bf u}_\mu^r
\right\}, \label{eq:quad_bp}
\end{eqnarray}

\noindent
where ${\bf a}^r_{\mu\rightarrow(\mu i)}$ denotes a vector excluding $a_{\mu r}$ from ${\bf a}_{\mu\rightarrow(\mu i)}=(a_{\mu 1\rightarrow(\mu i)}, ... , a_{\mu R\rightarrow(\mu i)})$, and $\Gamma_{{\bf a}^r_{\mu\rightarrow(\mu i)}}$ and $\Gamma_{{\bf b}^r_{\mu\rightarrow(\mu i)}}$ indicate, respectively, ${\rm diag}({\bf a}^r_{\mu\rightarrow(\mu i)}+\lambda{\bf 1})$ and ${\rm diag}({\bf b}^r_{\mu\rightarrow(\mu i)})$. Similarly for ${\bf b}^r_{\mu\rightarrow(\mu i)}$.

The minimization problem in (\ref{eq:quad_bp}) is solved as follows:

\begin{equation}
({\bf u}_\mu^r)^* = \left(\Gamma_{{\bf a}^r_{\mu\rightarrow(\mu i)}} + {\bf v}^r_i({\bf v}^r_i)^T\right)^{-1} \left\{{\bf b}_{\mu\rightarrow(\mu i)}^r - (y_{\mu i}-u_{\mu r}v_{ir}){\bf v}^r_i \right\}. \label{eq:quad_bp_solved}
\end{equation}

\noindent
Based on Sherman--Morrison formula, the inverse matrix in (\ref{eq:quad_bp_solved}) is re-expressed as follows:

\begin{equation}
(\Gamma_{{\bf a}^r_{\mu\rightarrow(\mu i)}} + {\bf v}^r_i({\bf v}^r_i)^T)^{-1} = \Gamma_{{\bf a}^r_{\mu\rightarrow(\mu i)}} - \frac{\Gamma_{{\bf a}^r_{\mu\rightarrow(\mu i)}}^{-1}{\bf v}^r_i({\bf v}^r_i)^T\Gamma_{{\bf a}^r_{\mu\rightarrow(\mu i)}}^{-1}}{1+({\bf v}^r_i)^T\Gamma_{{\bf a}^r_{\mu\rightarrow(\mu i)}}^{-1}{\bf v}^r_i}.\label{eq:sherma}
\end{equation}

\noindent
We insert (\ref{eq:quad_bp_solved}) and (\ref{eq:sherma}) into (\ref{eq:quad_bp}) to yield the following expression:

\begin{eqnarray}
\fl
\hat{f}_{(\mu i)\rightarrow\mu r}(u_{\mu r}) = \frac{1}{2}\frac{v_{ir}^2}{1+\chi_{(\mu i)}-\frac{v_{ir}^2}{a_{\mu r\rightarrow(\mu i)}+\lambda}}u_{\mu r}^2 - \frac{y_{(\mu i)} - \Delta_{(\mu i)} + u_{\mu r\rightarrow(\mu i)}v_{ir}}{1+\chi_{(\mu i)}-\frac{v_{ir}^2}{a_{\mu r\rightarrow(\mu i)}+\lambda}}v_{ir}u_{\mu r}
, \label{eq:fhat_cbmf}
\end{eqnarray}

\noindent
where $\chi_{(\mu i)}$, $\Delta_{(\mu i)}$ and $u_{\mu r\rightarrow(\mu i)}$ are defined as follows:

\begin{eqnarray}
\chi_{(\mu i)} &=& \sum_r \frac{v_{ir}^2}{a_{\mu r\rightarrow(\mu i)}+ \lambda}, \\
\Delta_{(\mu i)} &=& \sum_r u_{\mu r\rightarrow(\mu i)}v_{ir},
\label{eq:ab_mi}\\
u_{\mu r\rightarrow(\mu i)}&=&\frac{b_{\mu r\rightarrow(\mu i)}}{a_{\mu r\rightarrow(\mu i)}+ \lambda}.
\end{eqnarray}

\noindent
From the equations (\ref{eq:equad1}) and (\ref{eq:fhat_cbmf}), we obtain the following:

\begin{eqnarray}
\hat{a}_{(\mu i)\rightarrow\mu r}&=&\frac{v_{ir}^2}{1+\chi_{(\mu i)}-\frac{v_{ir}^2}{a_{\mu r\rightarrow(\mu i)}+\lambda}}, \\
\hat{b}_{(\mu i)\rightarrow\mu r}&=&\frac{y_{(\mu i)} - \Delta_{(\mu i)} + u_{\mu r\rightarrow(\mu i)}v_{ir}}{1+a_{(\mu i)}-\frac{v_{ir}^2}{\chi_{\mu r\rightarrow(\mu i)}+\lambda}}v_{ir} .
\end{eqnarray}

\noindent
Further, we insert (\ref{eq:equad1}) and (\ref{eq:equad3}) into (\ref{eq:bp3}) to yield the following expression:

\begin{eqnarray}
a_{\mu r\rightarrow(\mu i)}&=&a_{\mu r}-\hat{a}_{(\mu i)\rightarrow\mu r},\\
b_{\mu r\rightarrow(\mu i)}&=&b_{\mu r}-\hat{b}_{(\mu i)\rightarrow\mu r},
\end{eqnarray}

\noindent
where $a_{\mu r}$ and $b_{\mu r}$ are defined as follows:
\begin{eqnarray}
a_{\mu r} &=&\sum_{(\mu i)\in\partial\mu r} \hat{a}_{(\mu i)\rightarrow\mu r},\\
b_{\mu r} &=&\sum_{(\mu i)\in\partial\mu r} \hat{b}_{(\mu i)\rightarrow\mu r}.
\end{eqnarray}

\noindent
Finally, entries of the factorized matrices $u_{\mu r}^*$ are re-expressed from the equation (\ref{eq:fixedu}) as follows:
\begin{equation}
u_{\mu r}^*=\frac{b_{\mu r}}{a_{\mu r}+ \lambda}
\label{eq:u_mr_star}
\end{equation}

Similarly, we can re-express equations with respect to $\hat{c}_{(\mu i)\rightarrow ir},\hat{d}_{(\mu i)\rightarrow ir}$ and $c_{ir\rightarrow(\mu i)},d_{ir\rightarrow(\mu i)}$ based on (\ref{eq:bp2}),(\ref{eq:bp4}) and (\ref{eq:equad2}),(\ref{eq:equad4}).

In summary, the resulting equations are expressed as follows:

\noindent
\begin{itemize}
    \item Update equations for $U$:
\end{itemize}
\begin{eqnarray}
\chi^{t+1}_{(\mu i)} &=& \sum_r \frac{(v_{ir}^t)^2}{a^t_{\mu r\rightarrow(\mu i)}+ \lambda} \label{eq:cbmf1} \\
\Delta^{t+1}_{(\mu i)} &=& \sum_r u^t_{\mu r\rightarrow(\mu i)}v^t_{ir} \label{eq:cbmf2} \\
\hat{a}^{t+1}_{(\mu i)\rightarrow\mu r}&=&\frac{(v^t_{ir})^2}{1+\chi^t_{(\mu i)}-\frac{(v_{ir}^t)^2}{a^t_{\mu r\rightarrow(\mu i)}+\lambda}} \label{eq:a_hat} \\
\hat{b}^{t+1}_{(\mu i)\rightarrow\mu r}&=&\frac{y_{(\mu i)} - \Delta^t_{(\mu i)} + u^t_{\mu r\rightarrow(\mu i)}v^t_{ir}}{1+\chi^t_{(\mu i)}-\frac{(v_{ir}^t)^2}{a^t_{\mu r\rightarrow(\mu i)}+\lambda}}v^t_{ir}\label{eq:b_hat}\\
a^{t+1}_{\mu r} &=&\sum_{(\mu i)\in\partial\mu r} \hat{a}^t_{(\mu i)\rightarrow\mu r} \label{eq:a_mr_m}\\
b^{t+1}_{\mu r} &=&\sum_{(\mu i)\in\partial\mu r} \hat{b}^t_{(\mu i)\rightarrow\mu r} \label{eq:b_mr_m}\\
a^{t+1}_{\mu r\rightarrow(\mu i)}&=&a^t_{\mu r}-\hat{a}^t_{(\mu i)\rightarrow\mu r} \label{eq:a_mr}\\
b^{t+1}_{\mu r\rightarrow(\mu i)}&=&b^t_{\mu r}-\hat{b}^t_{(\mu i)\rightarrow\mu r} \label{eq:b_mr}\\
u^{t+1}_{\mu r\rightarrow(\mu i)} &=& \frac{b^t_{\mu r^\rightarrow(\mu i)}}{a^t_{\mu r\rightarrow(\mu i)}+ \lambda}\label{eq:cbmf_u_mi_}\\
u^{t+1}_{\mu r} &=& \frac{b^t_{\mu r}}{a^t_{\mu r}+ \lambda}\label{eq:cbmf_u}
\end{eqnarray}

\noindent
\begin{itemize}
    \item Update equations for $V$:
\end{itemize}
\begin{eqnarray}
\eta^{t+1}_{(\mu i)} &=& \sum_r \frac{(u_{\mu r}^{t+1})^2}{c^t_{ir\rightarrow(\mu i)}+ \lambda}\label{eq:c_mi} \\
\Theta^{t+1}_{(\mu i)} &=& \sum_r v^t_{\mu r\rightarrow(\mu i)}u^{t+1}_{\mu r} \\
\hat{c}^{t+1}_{(\mu i)\rightarrow ir}&=&\frac{(u_{\mu r}^{t+1})^2}{1+\eta^t_{(\mu i)}-\frac{(u_{\mu r}^{t+1})^2}{c^t_{ ir\rightarrow(\mu i)}+\lambda}} \\
\hat{d}^{t+1}_{(\mu i)\rightarrow ir}&=&\frac{y_{(\mu i)} - \Theta^t_{(\mu i)} + v^t_{ir\rightarrow(\mu i)}u^{t+1}_{\mu r}}{1+\eta^t_{(\mu i)}-\frac{(u_{\mu r}^{t+1})^2}{c^t_{ ir\rightarrow(\mu i)}+\lambda}}u^{t+1}_{\mu r}\\
c^{t+1}_{ir} &=&\sum_{(\mu i)\in\partial ir} \hat{c}^t_{(\mu i)\rightarrow ir}\\
d^{t+1}_{ir} &=&\sum_{(\mu i)\in\partial ir} \hat{d}^t_{(\mu i)\rightarrow ir}  \\
c^{t+1}_{ ir\rightarrow(\mu i)}&=&c^t_{ir}-\hat{c}^t_{(\mu i)\rightarrow ir}\\
d^{t+1}_{ ir\rightarrow(\mu i)}&=&d^t_{ir}-\hat{d}^t_{(\mu i)\rightarrow ir} \label{eq:d_ir}\\
v^{t+1}_{ir\rightarrow(\mu i)} &=& \frac{d^t_{ir\rightarrow(\mu i)}}{c^t_{ir\rightarrow(\mu i)}+ \lambda} \label{eq:cbmf_v_ir} \\
v^{t+1}_{ ir} &=& \frac{d^t_{ir}}{c^t_{ir}+ \lambda} \label{eq:cbmf20}
\end{eqnarray}

\noindent
Here, $t$ denotes the counter index for the update. It should be noted that in order to update variables for $V$ at time $t$, $u_{\mu r}^{t+1}$ is used instead of $u_{\mu r}^{t}$. 
We term the algorithm composed of (\ref{eq:cbmf1})-(\ref{eq:cbmf20}) as cavity-based matrix factorization (CBMF). 

The computational cost per update of each equation is $O(|\Omega |R)$ and the necessary memory cost corresponds to $O(|\Omega |R)$. The computational cost is competitive, and this is discussed later. Conversely, the necessary memory cost of CBMF exceeds those of ALS and SGD (Table \ref{tab:algorithms}). Although this is a disadvantage of CBMF, its necessary memory size is reduced to that of ALS and SGD by utilizing an approximation that is similar to that for deriving AMP from belief propagation \cite{kabashima2003cdma} as shown below.

\subsection{Derivation of the approximate algorithm\label{deriation_acbmf}}

CBMF entails $O(|\Omega |R)$ memory cost, and this is equivalent to the number of edges in the factor graph. 
When $R$ and $c$ are sufficiently large, the effect caused by omitting a variable node is expected to be negligible. Thus, the variables corresponding to the edges can be replaced by those corresponding to nodes.
The goal of this subsection involves deriving update equations with respect to the variables corresponding to the nodes. In the following, $R$ and $c$ are assumed as sufficiently large. 

The equation (\ref{eq:a_hat}) is approximately re-expressed as follows:

\begin{equation}
\hat{a}_{(\mu i)\rightarrow\mu r}
=\frac{v_{ir}^2}{1+\sum_s \frac{v_{is}^2}{a_{\mu s\rightarrow(\mu i)}+\lambda} -\frac{v_{ir}^2}{a_{\mu r\rightarrow(\mu i)}+\lambda}}
\simeq\frac{v_{ir}^2}{1+\chi_{(\mu i)}},
\label{eq:a_hat_mr_ap}
\end{equation}

\noindent
where $\chi_{(\mu i)}$ is also approximated by ignoring one of $c$ terms as follows:

\begin{equation}
    \chi_{(\mu i)}\simeq\sum_s \frac{v_{is}^2}{a_{\mu s}+\lambda}.
\label{eq:chi_mi}
\end{equation}

\noindent
Operating $\sum_{(\mu i)\in\partial\mu r}$ on both sides of (\ref{eq:a_hat_mr_ap}) yields

\begin{equation}
    a_{\mu r} =
    \sum_{(\mu i)\in\partial\mu r} 
    \frac{v_{ir}^2}{1+\chi_{(\mu i)}}.
\label{eq:app_a_mr}
\end{equation}

\noindent
Similarly, the equation (\ref{eq:b_hat}) is re-expressed as follows:

\begin{equation}
    \hat{b}_{(\mu i)\rightarrow\mu r}\simeq
    \left(
    \frac{y_{(\mu i)}-\sum_s u_{\mu s\rightarrow(\mu i)}v_{is}}{1+\chi_{(\mu i)}}
    + \frac{v_{ir}u_{\mu r\rightarrow(\mu i)}}{1+\chi_{(\mu i)}}
    \right)v_{ir},
\label{eq:ap_b_hat_mi_mr}
\end{equation}

\noindent
where $u_{\mu s\rightarrow(\mu i)}$ is also approximated by ignoring one of $R$ or $c$ terms as follows:

\begin{eqnarray}
    u_{\mu s\rightarrow(\mu i)}
    \simeq
    u_{\mu s} 
    - \phi_{(\mu i)}
    \frac{v_{is}}{a_{\mu s}+\lambda},
\label{eq:ap_u_ms_mi}
\end{eqnarray}

\noindent
where $\phi_{(\mu i)}$ is defined as follows:

\begin{eqnarray}
    \phi_{(\mu i)}
    &=&
    \frac{y_{(\mu i)}-\sum_{s} u_{\mu s\rightarrow(\mu i)}v_{is}}
    {1+\chi_{(\mu i)}} \\
    &\simeq&
    \frac{y_{(\mu i)}-\sum_{s} u_{\mu s}v_{is}+\phi_{(\mu i)}\chi_{(\mu i)}}
    {1+\chi_{(\mu i)}}.
\label{eq:ap_phi_mi}
\end{eqnarray}

\noindent
The second line is derived from (\ref{eq:chi_mi}) and (\ref{eq:ap_u_ms_mi}).
We insert (\ref{eq:ap_u_ms_mi}) into (\ref{eq:ap_b_hat_mi_mr}) and operate $\sum_{(\mu i)\in\mu r}$ on both sides to yield the following expression: 

\begin{equation}
    b_{\mu r}=
    \sum_{(\mu i)\in\mu r} \phi_{(\mu i)}v_{ir}
    + u_{\mu r}\sum_{(\mu i)\in\mu r} \frac{v_{ir}^2}{1+\chi_{(\mu i)}}.
\label{eq:ap_b_mr}
\end{equation}

Similarly, the update equations (\ref{eq:c_mi})-(\ref{eq:cbmf20}) are re-expressed by the same procedure. Finally, the approximate update equations are summarized as follows:

\begin{itemize}
    \item Update equations for $U$:
\end{itemize}
\begin{eqnarray}
    \chi^{t+1}_{(\mu i)}&=&\sum_s \frac{(v_{is}^t)^2}{a^t_{\mu s}+\lambda} \label{eq:approx_cmbf1} \\
    \phi^{t+1}_{(\mu i)}&=&
    \frac{y_{(\mu i)}-\sum_{s} u^t_{\mu s}v^t_{is}+\phi^t_{(\mu i)}\chi^t_{(\mu i)}}
    {1+\chi^t_{(\mu i)}}\\
    a^{t+1}_{\mu r} &=&
    \sum_{(\mu i)\in\partial\mu r} 
    \frac{(v_{ir}^t)^2}{1+\chi^t_{(\mu i)}}\\
    b^{t+1}_{\mu r}&=&
    \sum_{(\mu i)\in\mu r} \phi^t_{(\mu i)}v^t_{ir}
    + u^t_{\mu r}\sum_{(\mu i)\in\mu r} \frac{(v_{ir}^t)^2}{1+\chi^t_{(\mu i)}}\\
    u^{t+1}_{\mu r} &=& \frac{b^t_{\mu r}}{a^t_{\mu r}+ \lambda} \\
\end{eqnarray}

\begin{itemize}
    \item Update equations for $V$:
\end{itemize}
\begin{eqnarray}
    \eta^{t+1}_{(\mu i)}&=&\sum_s \frac{(u_{is}^{t+1})^2}{c^t_{i s}+\lambda} \\
    \psi^{t+1}_{(\mu i)}&=&
    \frac{y_{(\mu i)}-\sum_{s} u^{t+1}_{\mu s}v^t_{is}+\psi^t_{(\mu i)}\eta^t_{(\mu i)}}
    {1+\eta^t_{(\mu i)}}\\
    c^{t+1}_{i r} &=&
    \sum_{(\mu i)\in\partial ir} 
    \frac{(u_{\mu r}^{t+1})^{2}}{1+\eta^t_{(\mu i)}}\\
    d^{t+1}_{i r}&=&
    \sum_{(\mu i)\in ir} \psi^t_{(\mu i)}u^{t+1}_{\mu r}
    + v^t_{ir}\sum_{(\mu i)\in ir} \frac{(u_{\mu r}^{t+1})^{2}}{1+\eta^t_{(\mu i)}} \\
    v^{t+1}_{ir} &=& \frac{d^t_{ ir}}{c^t_{ ir}+ \lambda} \label{eq:approx_cmbf10}
\end{eqnarray}

\noindent
We term the algorithm composed of (\ref{eq:approx_cmbf1})-(\ref{eq:approx_cmbf10}) as the approximate cavity-based matrix factorization (ACBMF). The necessary memory cost to execute the algorithm is $O((N+M)R+|\Omega|)$, which is equivalent to the number of nodes in the factor graph. When compared to CBMF, ACBMF significantly reduces the required memory cost while the necessary computational cost is unchanged.

Additionally, one can illustrate that the fixed point of ACBMF is in agreement with that of ALS. Equation (\ref{eq:ap_phi_mi}) is solved with respect to $\phi_{(\mu i)}$, and we obtain the following expression:

\begin{equation}
    \phi_{(\mu i)}=y_{(\mu i)}-\sum_s u_{\mu s}v_{is}.
\label{eq:phi_mi_fix}
\end{equation}

\noindent
We insert (\ref{eq:app_a_mr}) and (\ref{eq:phi_mi_fix}) into (\ref{eq:ap_b_mr}) to yield the following expression:

\begin{equation}
    b_{\mu r}= \sum_{(\mu i)\in\mu r} \left(y_{(\mu i)}-\sum_s u_{\mu s}v_{is}\right)v_{ir}+u_{\mu r}a_{\mu r}.
\end{equation}

\noindent
From the equations (\ref{eq:u_mr_star}), we obtain the following expression:

\begin{equation}
    \lambda{\bf u}_r=\sum_{(\mu i)\in\mu r}
    \left(
        y_{(\mu i)}-{\bf u}^T_r{\bf v}_i
    \right){\bf v}_i.
\label{eq:lambda_u_r}
\end{equation}

\noindent
We solve (\ref{eq:lambda_u_r}) with respect to ${\bf u}_r$ to yield the following expression:

\begin{equation}
    {\bf u}^*_r = \left(
        \sum_{(\mu i)\in\mu r}{\bf v}_i{\bf v}^T_i+\lambda{\bf I}_R
    \right)^{-1}
    \left(\sum_{(\mu i)\in\mu r} y_{(\mu i)}{\bf v}_i\right),
\label{eq:als_inv2}
\end{equation}

\noindent
and this is equivalent to (\ref{eq:als_inv}). Similarly for ${\bf v}^*_r$.

In contrast to ALS, ACBMF does not completely optimize $U$ ($V$) for a given $V$ ($U$) in each step, and thus the necessary computation is reduced. Evidently, this may decrease the convergence speed. However, the complete optimization for it does not necessarily bring $U$ ($V$) to a better state when $V$ ($U$) is far from the convergent solution. Therefore, it is not advised to expend significant computational cost on this. Additionally, the optimization in each step tends to strengthen time correlations of the variables, and this may make the cavity treatment inappropriate. Actually, the results of experiments shown below indicate that this concern is the case.

\subsection{Comparison with ALS and SGD}
We briefly compare (A)CBMF with ALS and SGD. ALS and SGD are algorithms that attempt to iteratively minimize the multivariate objective function (\ref{eq:nonconvex_relaxed}). Although their working principle is natural, the performance of these algorithms can be negatively affected by the self-feedback effect caused by cycles from the graph. Conversely, (A)CBMF reduces such effect by introducing the seemingly artificial cavity functions, and this may lead to the performance improvement. In a manner similar to ALS, (A)CBMF can also be easily parallelized, and is free from learning parameters unlike SGD.

The computational and memory costs of the four algorithms are summarized in Table \ref{tab:algorithms}. The computational cost is defined as that necessary to update all variables at least once. Given this definition, SGD only updates the variables based on the gradients although the computational cost of SGD appears the lowest. Conversely, CBMF and ALS update them with closed forms, and thus it is expected that their convergence speeds can increase. A comparison of (A)CBMF and ALS indicates that the computational cost of the former is lower. Conversely, the memory cost of CBMF is the highest while that of ACBMF is identical to that of ALS and SGD.

\begin{table}[htb]
\begin{center}
\scalebox{0.81}{
  \begin{tabular}{|l|c|c|c|c|} \hline
     & CBMF & ACBMF & ALS & SGD \\ \hline \hline
    Computational costs & $O(|\Omega|R)$ &　$O(|\Omega|R)$ & $O((|\Omega|R^2+(N+M)R^3))$&   $O((N+M)R)$   \\ \hline
    Memory costs & $O(|\Omega|R)$ & $O((N+M)R+|\Omega|)$ & $O((N+M)R+|\Omega|)$&   $O((N+M)R+|\Omega|)$   \\ \hline
  \end{tabular}
  }
  \caption{Comparison of computational costs to update all variables at least once. Specifically, $|\Omega|$ denotes the number of observed entries, and this is assumed to exceed or be equal to the number of variables to be determined $(N+M)R$.}
  \label{tab:algorithms}
\end{center}
\end{table}

\section{Numerical Experiments}

\subsection{Synthetic Data Analysis}

In order to systematically compare the performance of the four algorithms, namely ALS, SGD, and (A)CBMF (C++ implementation is available at \cite{cbmf_url}), we performed extensive numerical experiments using synthetic datasets. A dataset for the experiment was prepared as follows: The original matrix $Y^0\in \mathbb{R}^{N\times M}$ is provided from $U^0\in\mathbb{R}^{N\times R},V^0\in\mathbb{R}^{M\times R}$, and $Z\in\mathbb{R}^{N\times M}$ as $Y^0=U^0(V^0)^T+Z$, where entries of $U^0$ and $V^0$ are independently sampled from the standard Gaussian distribution while those of $Z$ are independently and identically distributed based on a Gaussian of zero mean and variance 0.09. We randomly select ``observed entries'' out of $Y^0$ with probability of $c/N$ where $c\sim O(1)$ denotes the average number of the observed entries per column. The collection of the observed entries constitutes the observed matrix $Y$. We assume that true rank $R$ is known in advance.

We evaluate the performance of the algorithms via the relative root mean square error (rRMSE) and the reconstruction rate. Given the effect of the noise $Z$, it is impossible to perfectly reconstruct $Y^0$ in the current setting. Therefore, we consider estimated factorized matrices $U$ and $V$ as successful if $\sqrt{\sum_{\mu i}(y_{\mu i}^0 - u_{\mu r}v_{ir})^2}/\sqrt{\sum_{\mu i}(y_{\mu i}^0)^2}\leq0.15$ holds. The convergence of the three algorithms is not guaranteed, and thus we attempt ten random initial conditions for each sample and algorithm and counted a ``success'' if at least one of the ten initial conditions leads to the successful reconstruction. Additionally, rRMSE is evaluated via the mean of the minimum value of $\sqrt{\sum_{\mu i}(y_{\mu i}^0 - u_{\mu r}v_{ir})^2}/\sqrt{\sum_{\mu i}(y_{\mu i}^0)^2}$ out of the ten initial conditions over 50 samples. Conversely, the reconstruction rate denotes the fraction of the reconstruction success over the 50 samples.

\begin{figure}[h]
\centering
\begin{tabular}{c}
\hspace{-20pt}
\subfigure[]{
\includegraphics[clip, width=8.6cm]{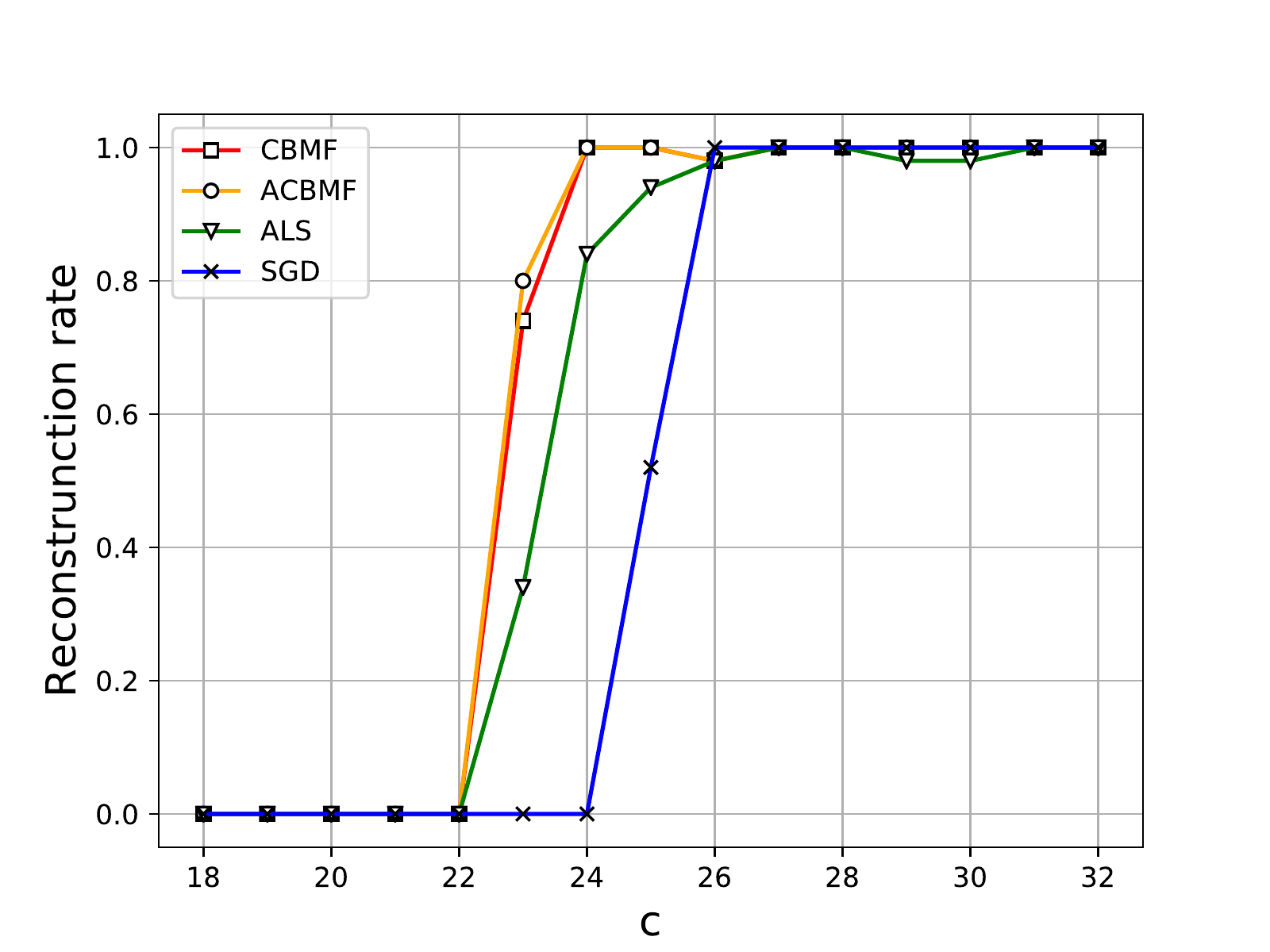}
\label{fig:synthetic_recon}}
\hspace{-30pt}
\subfigure[]{
\includegraphics[clip, width=8.6cm]{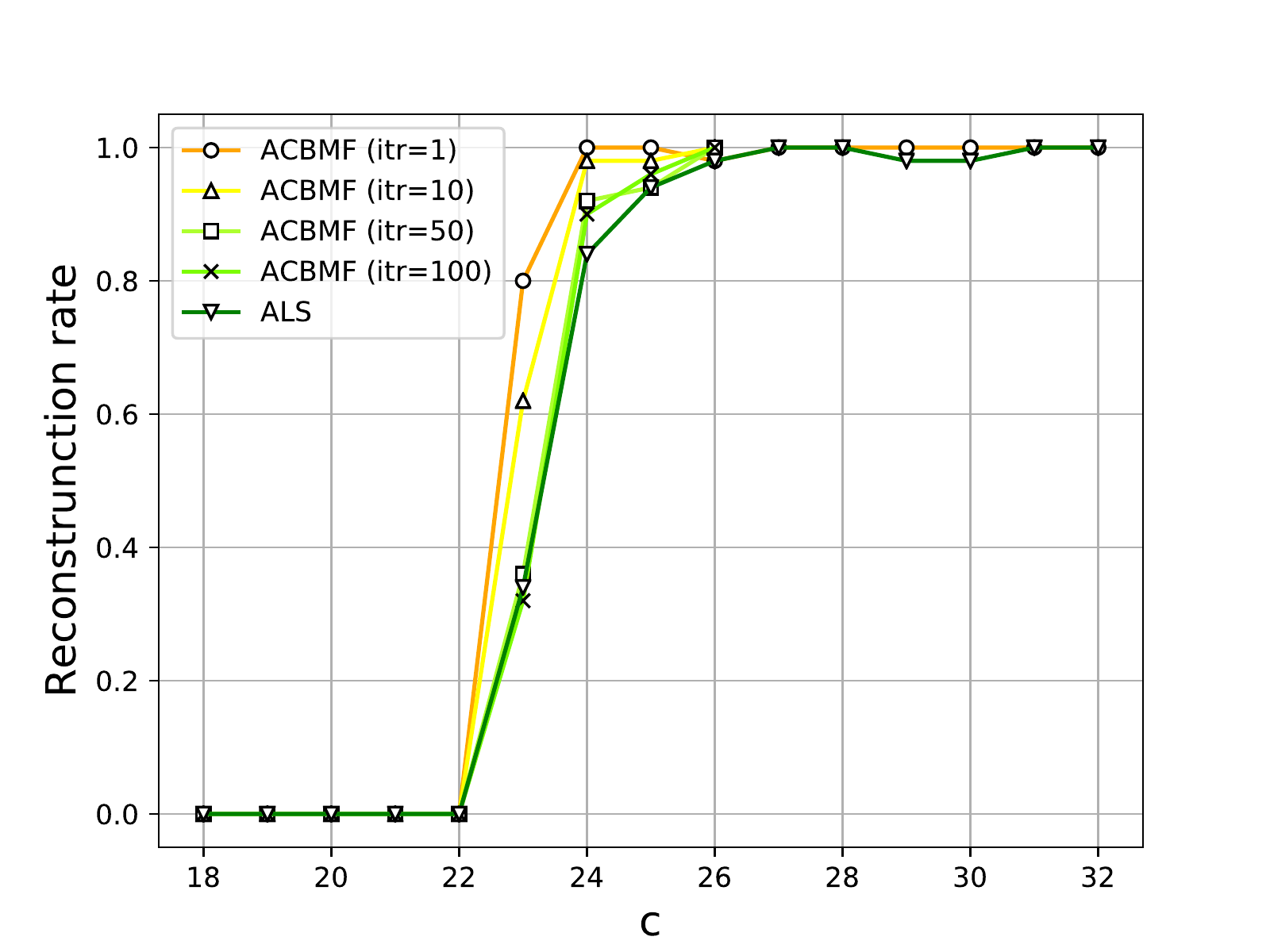}
\label{fig:synthetic_recon_local_iteration}}
\end{tabular}
\caption{Reconstruction rate as a function of $c$ for matrices with rank $R=10$, and system size $N=500, M=1000$. For each $c$, the rate was evaluated from 50 experiments where $\lambda=10^{-2}$ was used. Reconstruction is considered as successful when $\sqrt{\sum_{\mu i}(y_{\mu i}^0 - u_{\mu r}v_{ir})^2}/\sqrt{\sum_{\mu i}(y_{\mu i}^0)^2}\le0.15$ for at least once in ten trials.
(a) Comparison between (A)CBMF, ALS and SGD. (b) Results for ACMBF when the number of iterations for each step increases.}
\label{fig:synthetic_recon_both}
\end{figure}

Figure \ref{fig:synthetic_recon} plots the experimental results as function of the average number $c$ of observations per column for $R=10$. The figure indicates that (A)CBMF outperforms the other algorithms. It should be noted that (A)CBMF exhibits a better reconstruction rate up to a smaller value of $c$ than ALS while they are theoretically guaranteed to share the same fixed point. We speculate that this is because (A)CBMF weakens the self-feedback effect via the cavity treatment and by not performing optimization in each step. In order to verify the validity of this speculation, we examine the manner in which the reconstruction rate changes when the number of iterations of ACBMF for each step increases, and this is plotted in figure \ref{fig:synthetic_recon_local_iteration}. When the iteration is repeated until convergence in each step, $U$ ($V$) is optimized for a given $V$ ($U$). This implies that the performance would become {\em worse} when the number of the iterations increases by spending more computational cost. The figure shows that this is actually the case and supports our speculation.

\begin{figure}[h]
\centering
\begin{tabular}{c}
\hspace{-20pt}
\subfigure[]{
\includegraphics[clip, width=8.6cm]{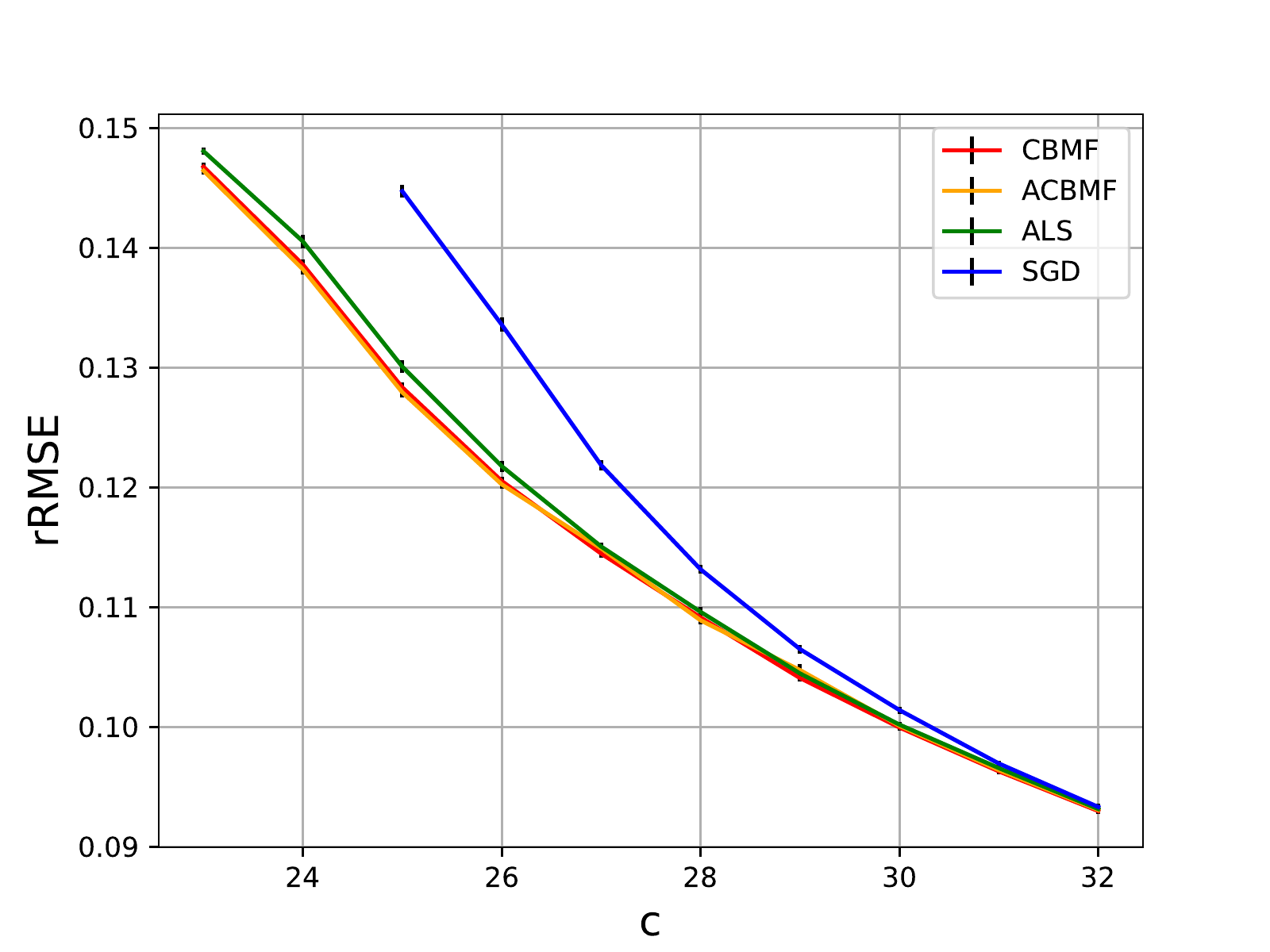}
\label{fig:synthetic_rrmse}}
\hspace{-30pt}
\subfigure[]{
\includegraphics[clip, width=8.6cm]{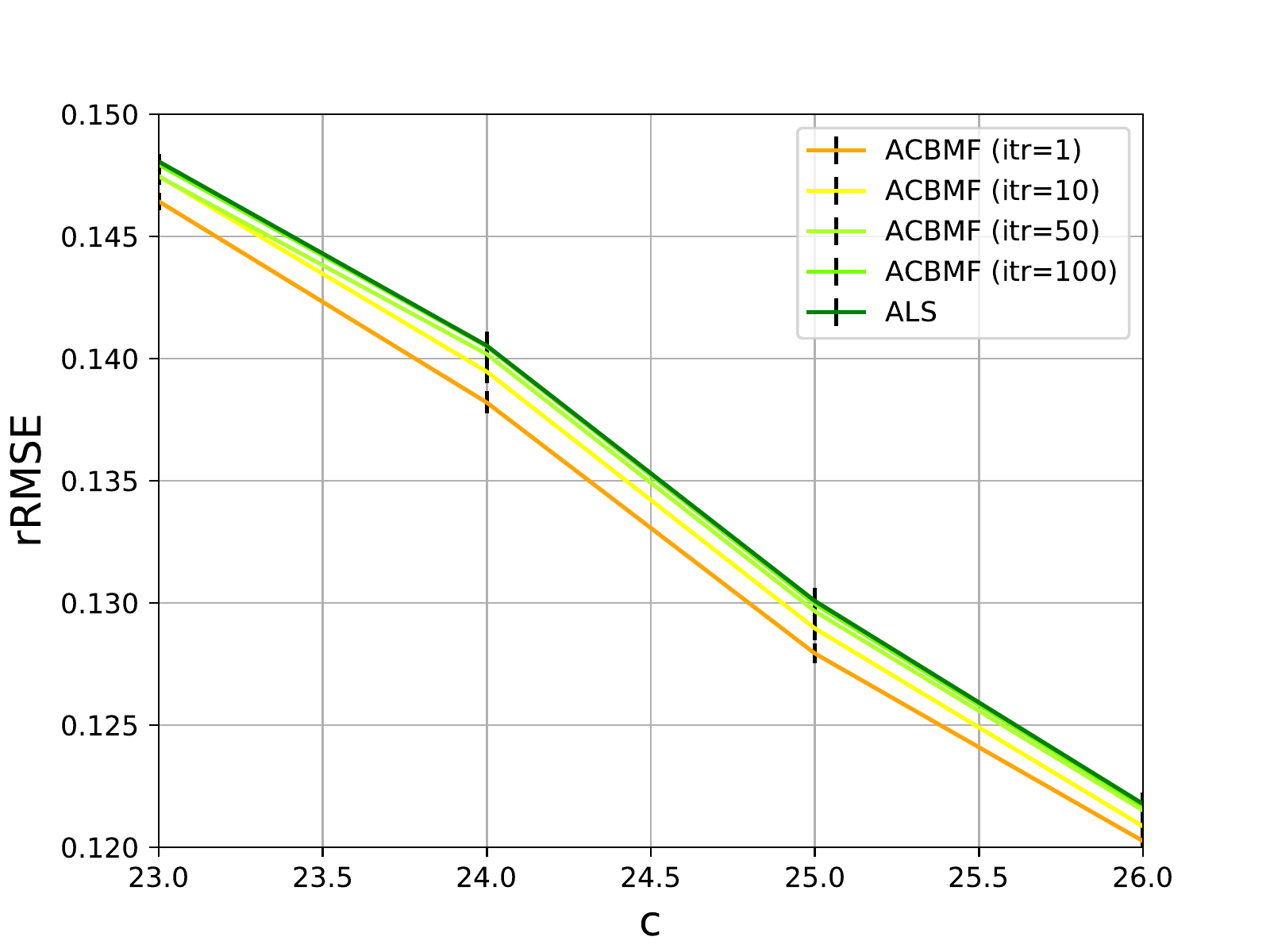}
\label{fig:synthetic_rrmse_local_iteration}}
\end{tabular}
\caption{Relative root mean square error (rRMSE) of reconstructed samples as a function of $c$. Experimental conditions are identical to those in Figure \ref{fig:synthetic_recon_both}. (a) Comparison between (A)CBMF, ALS and SGD. (b) Results for ACBMF when the number of iterations for each step increases.}
\label{fig:synthetic_rrmse_both}
\end{figure}

Figure \ref{fig:synthetic_rrmse} shows the results for rRMSE. The performance of SGD is significantly worse when compared to that of (A)CBMF and ALS. This is potentially because the scheduling of the learning rate used in the SGD experiments is not optimally tuned. The default scheduling that is provided in a code distribution \cite{sgd_url} leads to a terrible result, and thus we select a better scheduling although it is non-trivial to determine the optimal one. Conversely, (A)CBMF and ALS are free from such issues as they involve no scheduling of parameters. (A)CBMF exhibits slightly better performance when compared to that ALS. Similarly, for reconstruction rate, the performance of ACBMF approaches that of ALS when the number of iterations per update increases (figure\ref{fig:synthetic_rrmse_local_iteration}).

\subsection{Real Data Analysis}

We also examined the usefulness of the proposed algorithm via application to three benchmark datasets of recommender systems, namely MovieLens 1M, 10M, and 20M \cite{harper2016movielens}. Specifically, the 1M dataset is composed of rating values $s$ from 1 to 5 with step 1, and 10M and 20M are from 0.5 to 5 with step 0.5. The higher values correspond to higher evaluations for movies or music provided by users. Details of the datasets are summarized in Table\ref{tab:datasets}.

The performance of each algorithm for the matrix is evaluated as follows: We randomly split the matrix entries into 10 groups, matrix factorization is performed by using data of 9-of-the-10 groups, and the performance of the obtained factorization is measured by using data of the remaining group. We employ root mean square error (RMSE) as a performance measure, and it is averaged over 50 samples of the experiment. In all the experiments, we set $R=10$.

\begin{figure}[htbp]
  \begin{center}
    \includegraphics[width=13cm]{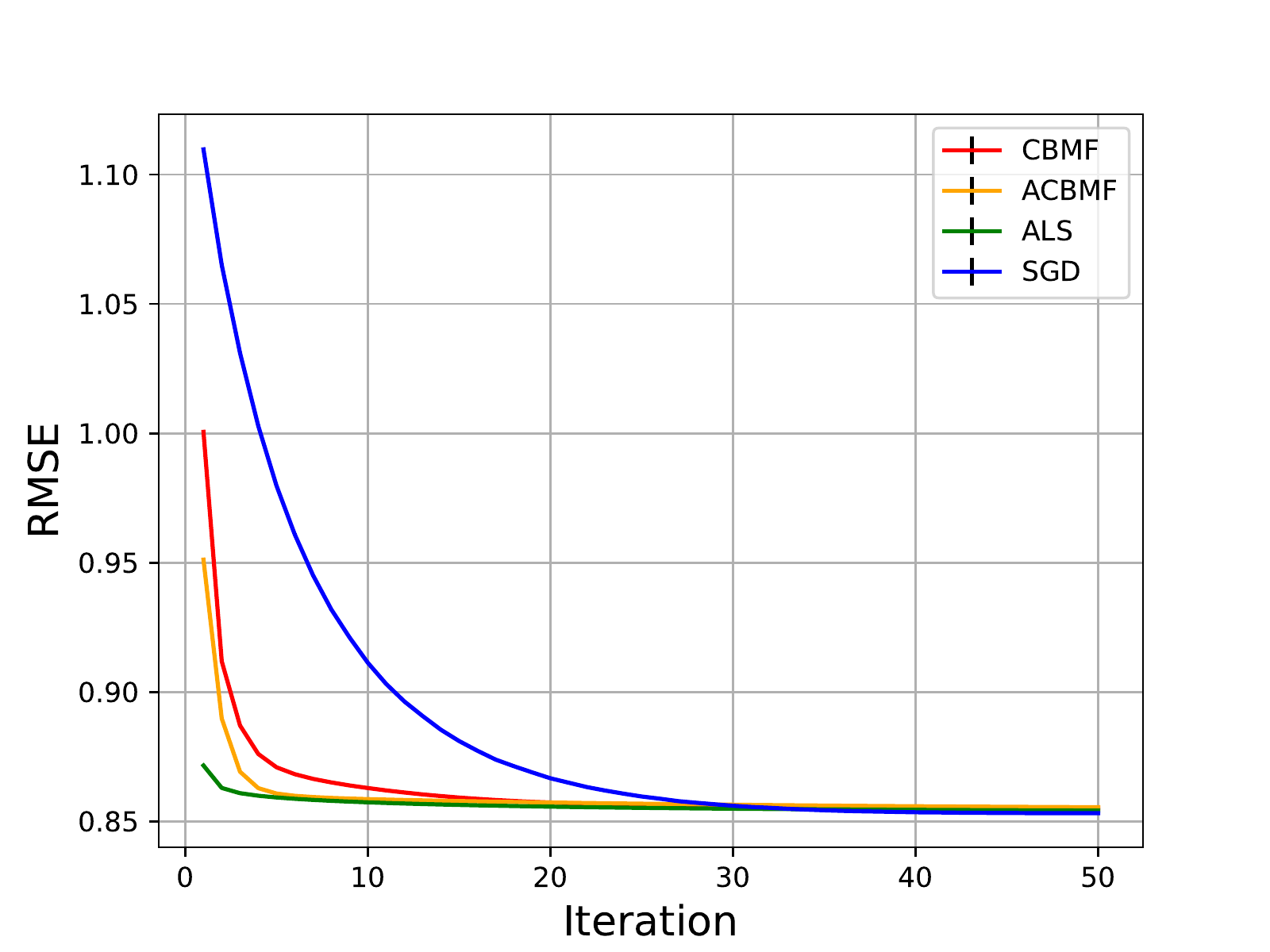}
    \caption{Results for MovieLens 1M dataset, RMSE is plotted versus iteration. In the experiments, we set $R=10$, and the regularization parameter $\lambda$ is fixed as 3. The figure compares the result of the four algorithms, (A)CBMF, ALS, and SGD.}
    \label{fig:iteration_rmse_1m}
  \end{center}
\end{figure}

\begin{figure}[htbp]
  \begin{center}
    \includegraphics[width=13cm]{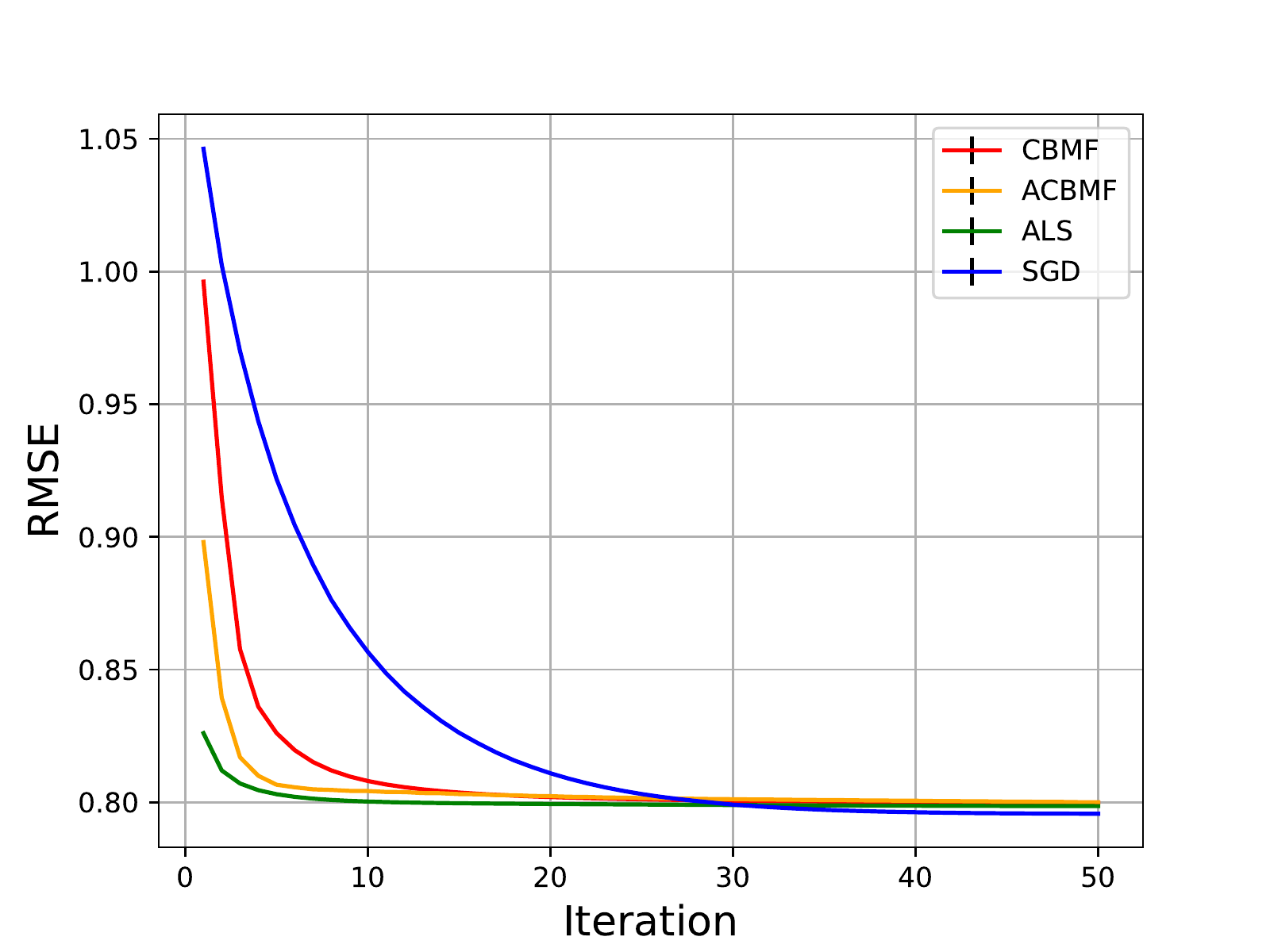}
    \caption{Results for MovieLens 10M dataset. Experimental conditions are identical to those in Figure \ref{fig:iteration_rmse_1m}. }
    \label{fig:iteration_rmse_10m}
  \end{center}
\end{figure}

\begin{figure}[htbp]
  \begin{center}
    \includegraphics[width=13cm]{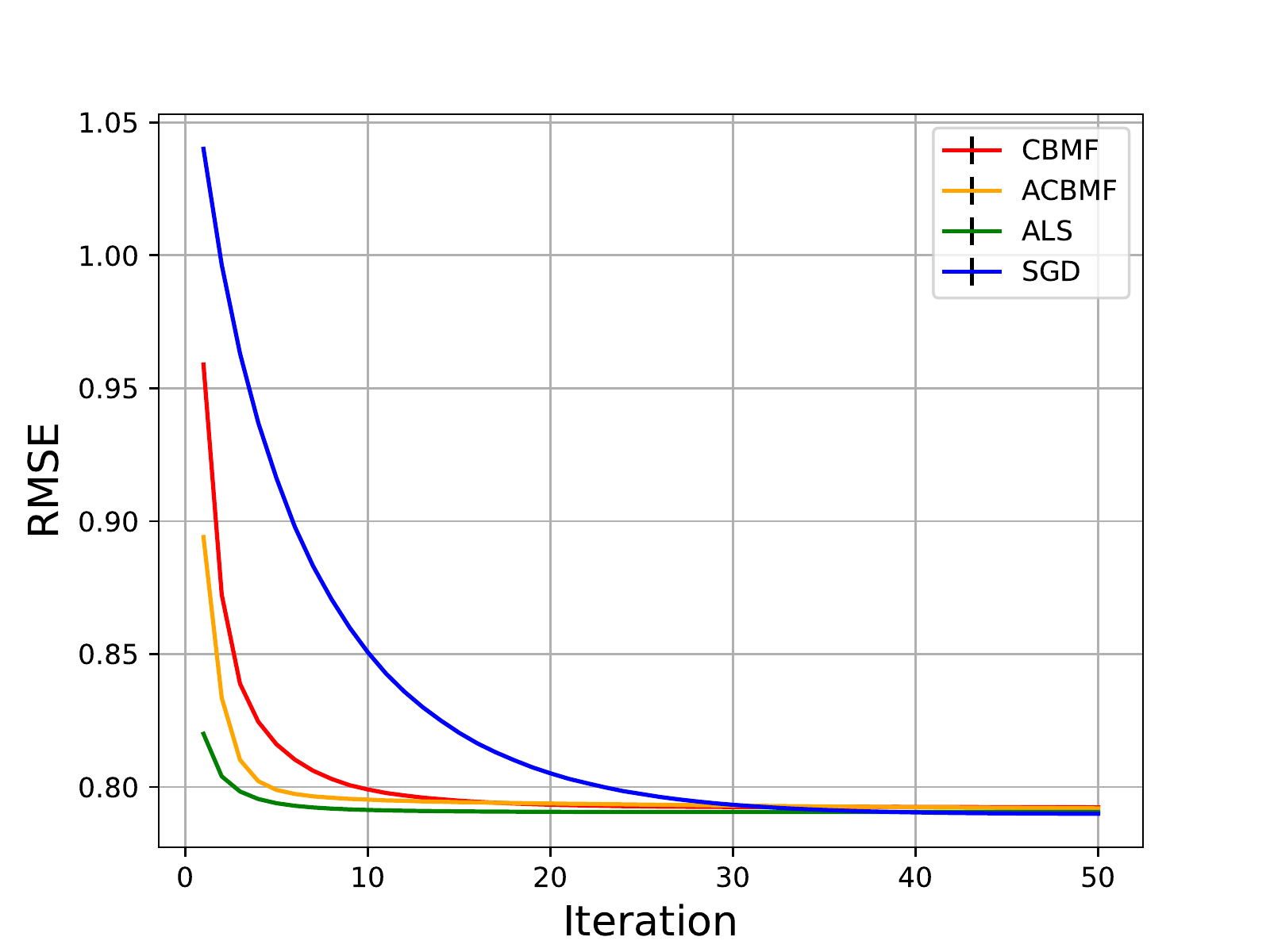}
    \caption{Results for MovieLens 20M dataset. Experimental conditions are identical to those in Figure \ref{fig:iteration_rmse_1m}. }
    \label{fig:iteration_rmse_20m}
  \end{center}
\end{figure}

Figures \ref{fig:iteration_rmse_1m}-\ref{fig:iteration_rmse_20m} show the performance measure of (A)CBMF, ALS and SGD evaluated for the three datasets. The figures represent RMSE relative to the number of iterations. The figures indicate that all the algorithms finally achieve similar performance although the number of iterations necessary for convergence is minimized for ALS. However, it should be noted that the ALS requires a significantly higher computational cost than (A)CBMF and SGD per iteration (Table\ref{tab:algorithms}). Thus, (A)CBMF converges faster than the other algorithms in terms of actual time when $R$ is relatively large.

\section{Summary}
In summary, we developed matrix factorization algorithms that are abbreviated as CBMF and ACBMF based on the cavity method. In terms of computational cost, CBMF is competitive with SGD because CBMF updates variables in closed forms (which generally reduces the number of iterations necessary for convergence) although a comparison of the necessary computational cost to update all variables at least once indicates that the computational cost of SGD is the smallest of the three. In a manner similar to CBMF, ALS updates variables in closed form although its computational cost exceeds that of CBMF because ALS requires the matrix inversion operation, which CBMF does not require. Conversely, in terms of the memory cost, CBMF requires more capacity than the others, and thus we developed ACBMF by utilizing an approximation that is similar to that for deriving AMP from belief propagation. The necessary memory cost of ACBMF is identical to that of SGD and ALS.

Experiments involving synthetic data indicated that (A)CBMF exhibits better performance without the necessity of parameter tuning when observed entries are not sufficiently large. The superiority of the performance presumably stems from the reduction of self-feedback effects via the introduction of cavity treatment and avoidance of the complete optimization in each update.
Experiments using real world dataset indicated that all algorithms achieved similar performance although (A)CBMF converges faster than the other two in actual time when rank $R$ is relatively large.

Future work includes generalization of CBMF to matrix factorization problems with additional constraints such as non-negative matrix factorization \cite{lee1999learning}.

\section*{Acknowledgements}
Useful discussion with Tomoyuki Obuchi is acknowledged. This study was partially supported by KAKENHI No.17H00764.

\clearpage
\appendix

\section{Benchmark datasets}

We performed numerical experiments on three different benchmark datasets as follows: the MovieLens 1M, 10M, and 20M datasets (https://movielens.org/). The characteristics of each dataset is represented in Table \ref{tab:datasets}.

\begin{table}[htb]
\begin{center}
  \begin{tabular}{|l|c|r|r|rl} \hline
    Dataset & Rating set & \#Users & \#Items & \#Ratings \\ \hline \hline
    MovieLens 1M & \{1,2,3,4,5\} & 6,040 & 3,900 & 1,000,209 \\
    MovieLens 10M & \{0.5,1,1.5,2,2.5,3,3.5,4,4.5,5\} & 10,681 & 71,567 & 10,000,054 \\
    MovieLens 20M & \{0.5,1,1.5,2,2.5,3,3.5,4,4.5,5\} & 138,493 & 27,278 & 20,000,263 \\ \hline
  \end{tabular}
  \caption{The details of the datasets used in this study. MovieLens is a dataset that consists of the ratings for movies from users who watched the movies, and the ratings of 1M dataset takes an integer value from 1 to 5 and those of 10M and 20M datasets take a value from 0.5 to 5 with step 0.5. When a user likes a movie very much, he or she rates  the movie as 5. \#Users and \#Items correspond to the row and column sizes of the observed matrix, respectively, and \#Ratings denotes the number of observations.}
  \label{tab:datasets}
\end{center}
\end{table}


\bibliographystyle{unsrt}
\bibliography{mybib}

\end{document}